\documentclass[12pt]{amsart}
\usepackage{amsgen,amsmath,amstext,amsbsy,amsopn,amsfonts,amssymb}
\usepackage{amsmath,amssymb,amsfonts,amsthm,enumerate}
\usepackage{mathrsfs}

\newtheorem{theorem}{Theorem}

\newtheorem{corollary}{Corollary}
\newtheorem{proposition}{Proposition}
\newtheorem{obs}{Observation}

 \newtheorem{defi}{Definition}
\newenvironment{definition}{\begin{defi}\rm}{\end{defi}}
\newtheorem{exa}{Example}

\newtheorem{rem}{Remark}
\newenvironment{remark}{\begin{rem}\rm}{\end{rem}}
\newtheorem{rems}{Remarks}

\newtheorem{ack}{Acknowlegment}

\DeclareMathOperator{\Span}{span}

\def\proof{\noindent\textbf{Proof}\quad}
\def\bsq{\blacksquare\medskip}
\def\n{\noindent}
\def\H{\mathcal H}
\def\N{\mathcal N}

\def\K{\mathcal K}
\def\P{\mathcal P}
\def\M{\mathcal M}

\def\B{\mathcal B}

\def\NN{{\mathbb N}}
\def\ZZ{{\mathbb Z}}

\def\FF{\mathbb F}
\def\RR+{{\mathbb R}^*}

\def\Q_p{{\mathbb Q}_p}

\def\Ga{\Gamma}

\def\al{\alpha}

\def\lb{\label}
\def\tout{\qquad\text{for all}\quad}
\def\n{\noindent}
\def\har{{\mathrm Har}_\mu}
\def\Har{{\mathrm Har}}
\def\Hom{{\mathrm Hom}}

\begin{document}

\title[Harmonic cocycles and  irreducible isometric actions]{Harmonic cocycles, von Neumann algebras, and  irreducible affine isometric actions}
\author{Bachir Bekka}
\address{Bachir Bekka \\ IRMAR \\ UMR-CNRS 6625 Universit\'e de  Rennes 1\\
Campus Beaulieu\\ F-35042  Rennes Cedex\\
 France}
\email{bachir.bekka@univ-rennes1.fr}

\thanks{The  author acknowledges the partial support of the  French Agence Nationale de la Recherche (ANR)
through the projects Labex Lebesgue (ANR-11-LABX-0020-01) and   GAMME  (ANR-14-CE25-0004). }

\begin{abstract}
Let $G$ be a  compactly generated  locally compact group and $(\pi, \H)$ a unitary representation of $G.$
The $1$-cocycles  with coefficients in $\pi$ which are harmonic (with respect to a suitable probability measure on $G$)
represent classes in the first reduced cohomology $\bar{H}^1(G,\pi).$
We show that harmonic $1$-cocycles 
 are characterized  inside their reduced cohomology class by the fact that they span a minimal closed 
subspace of $\H.$ In particular, the affine isometric action  given by a harmonic cocycle $b$ is
irreducible (in the sense that $\H$ contains no non-empty, proper closed invariant affine subspace)  if the linear span of $b(G)$ is dense in $\H.$ The converse statement is true,
if $\pi$ moreover has no almost invariant vectors.
Our approach exploits the natural structure of 
the space of harmonic $1$-cocycles with coefficients in $\pi$ 
as a Hilbert  module over  the von Neumann algebra $\pi(G)',$ which is the commutant of $\pi(G)$.
Using operator  algebras techniques, such as the von   Neumann dimension, 
we give a necessary and sufficient condition for a factorial representation $\pi$
 without almost invariant vectors to admit an irreducible affine  action  with  $\pi$ as linear part.
\end{abstract}
\maketitle

\section{Introduction}
\lb{S:Intro}
Let $G$ be a locally compact group and $(\pi, \H)$ a continuous unitary (or orthogonal) representation of
$G$ on a complex (or real) Hilbert space $\H.$ 
Recall that a \emph{$1$-cocycle} with coefficients  in $\pi$
is a continuous map $b: G\to \H$ such that 
$b(gh)= b(g)+ \pi(g)b(h)$ for all $g,h\in \H$ 
and that a $1$-cocycle is a \emph{coboundary}  if it is of the form 
$\partial_v$ for some $v\in \H,$ where $\partial_v(g)=\pi(g)v-v$ for $g\in G.$
The space  $Z^1(G,\pi)$  of $1$-cocycles with coefficients  in $\pi$ is a 
vector space containing the space $B^1(G,\pi)$ of coboundaries as linear subspace. 
The \emph{$1$-cohomology} $H^1(G,\pi)$ is the quotient $Z^1(G,\pi)/B^1(G,\pi)$.

The space $B^1(G,\pi)$ is not necessarily closed in $Z^1(G,\pi)$ (see Proposition~\ref{Pro-Gui})
and  the  \emph{reduced $1$-cohomology} with  coefficients in $\pi$ is defined as 
  $\overline{H}^1(G,\pi)= Z^1(G,\pi)/\overline{B^1(G,\pi)}.$
  
Assume now that $G$ is compactly generated,
 that is,  $G=\cup_{n\in \ZZ} Q^n$  for a compact
 subset $Q$, which we can assume to be a neighbourhood of the identity $e\in G$ 
 and to be symmetric ($Q^{-1}=Q$). 
 
  Harmonic $1$-cocycles in $Z^1(G,\pi),$  with respect to an appropriate probability measure on 
$G$, form a set of representatives for the classes in 
the reduced cohomology $\overline{H}^1(G,\pi),$ as we will shortly explain.
Such  cocycles appear in \cite{BV}   in the case where $\pi$ is the regular representation of a discrete group $G,$ 
in relation with the first $\ell^2$-Betti number of $G$; 
they play an important role in  Ozawa's recent  proof of Gromov's polynomial growth theorem
(\cite{Ozawa}) as well as  in the work \cite{Anna-Ozawa} and \cite{Gournay-Jolissaint}.

Harmonic $1$-cocycles were  implicitly introduced in \cite[Theorem 2]{Gui2};
it was  observed there that  $Z^1(G,\pi)$ can be identified with a  closed subspace of the Hilbert space $L^2(Q, \H,m_G)$, where
  $m_G$ is a (left) Haar measure on $G$ and so  $\overline{H}^1(G,\pi)$ corresponds to the orthogonal complement $B^1(G,\pi)^\perp$  of  $B^1(G,\pi)$ in $Z^1(G,\pi).$
  Following \cite{Anna-Ozawa}, we prefer to embed $Z^1(G,\pi)$ in a 
  more general  Hilbert space,  defined by  a class of appropriate probability measures similar to those appearing
  there.   For this, we consider the  word length on $G$ associated to $Q$, that is, the map $g\mapsto |g|_Q$, where
  $$|g|_Q= \min\{ n\in \NN\, : g\in Q^n\}.$$

 \begin{definition}
 \label{Def-AdmissibleMeasure}
 A probability measure $\mu$ on $G$ is 
 \emph{cohomologically adapted} (or, more precisely, $1$-cohomologically adapted) if it has  the following properties:
   \begin{itemize}
 \item $\mu$ is symmetric;
  \item $\mu$ is absolutely continuous with respect to the Haar measure $m_G;$
  \item  $\mu$ is adapted:  the support of $\mu$    is a generating set for $G;$
  \item $\mu$ has a second moment: $\int_{G} |x|_Q^2 d\mu(x) <\infty.$
    \end{itemize} 
  \end{definition}
     Observe that the class of cohomologically adapted measures is independent of  the generating compact set $Q$, 
since the length functions associated to two compact generating sets are bi-Lipschitz equivalent.   
    
 We consider the Hilbert space $ L^2(G,\H, \mu)$ of measurable square-integrable maps  $F: G\to \H.$ 
  Then $ Z^1(G,\pi)$ is a subset of  $L^2(G,  \H, \mu)$ (see Section~\ref{S:Cocycle-VN}). 
 Moreover, the linear operator $$\partial: \H \to  Z^1(G,\pi), v\mapsto \partial_v$$ is bounded, has 
 $B^1(G,\pi)$ as range, and it is straightforward to check that its adjoint 
 is $-\dfrac{1}{2} M_\mu,$ where
 $$M_\mu: Z^1(G,\pi) \to \H,\   b\mapsto \int_G b(x) d\mu(x).
 $$
 So, the orthogonal complement $B^1(G,\pi)^\perp$ of $B^1(G,\pi)$ in $Z^1(G,\pi)$ can be identified with 
 the space of harmonic cocycles in the sense of the following definition.
 In particular, the reduced cohomology $\overline{H}^1(G,\pi)$ can be identified with  $\har(G,\pi).$
 \begin{definition}
 \lb{HarmCocycle}
 A cocycle $b\in Z^1(G,\pi)$ is $\mu$-\emph{harmonic} 
 if $M_\mu(b)=0,$ that is, $\int_G b(x) d\mu(x) =0.$
 We denote by $\har(G,\pi)$ the space of $\mu$-harmonic cocyles in $Z^1(G,\pi)$ and by 
 $$P_{\Har}: L^2(G,\H,\mu) \to \har(G,\pi)$$
  the orthogonal projection on $\har(G,\pi).$
 \end{definition} 
Observe that, by the cocycle relation, $b\in Z^1(G,\pi)$ is $\mu$-harmonic if and only if it has the mean value
 property
 $$
 b(g)=\int_{G} b(gx) d\mu(x) \tout g\in G.
 $$

 In our opinion, the Hilbert space structure of $\overline{H}^1(G,\pi)$ given by its realization as a space of harmonic cocycles, together with its module structure over the von Neumann algebra
 $\pi(G)'$ (see below),  deserves more attention than it has received so far in the literature.
Our aim in this paper is to use this structure in relation with a natural notion of irreducibility for affine isometric actions
(see Definition~\ref{Def-IAA}).

Our first result shows that harmonic $1$-cocycles $b$ are characterized by 
 a remarkable minimality  property of the space $\overline{\Span(b(G))},$  
 the closure of the linear span of $b(G).$
  
  \begin{theorem}
  \lb{Theo1}
 Let $G$ be a  compactly generated group.
  Let $(\pi, \H)$ be an orthogonal or unitary representation of $G$ 
  and $\mu$ a cohomologically adapted  probability measure   on $G$.
 Let $b\in \har(G, \pi)$ be a $\mu$-harmonic cocycle.
 We have 
 $$\overline {\Span (b(G))}=\bigcap_{b'} \overline {\Span (b'(G))},$$
 where $b'$ runs over the $1$-cocycles in the cohomology class of $b$ in $\overline{H}^1(G,\pi).$
  \end{theorem}
   In particular,  Theorem~\ref{Theo1}  shows that,  for a $\mu$-harmonic cocycle $b,$ the 
  closed linear subspace spanned by  $b(G)$ only  depends 
 on the reduced cohomology class of $b$  and not on the choice of $\mu.$

 Recall that, given a  cocycle $b\in Z^1(G,\pi),$ a continuous   action $\alpha_{\pi, b}$ of $G$  on 
 ${\mathcal H}$  by affine isometries is defined  by the formula 
$$\alpha_{\pi,b}(g)v=\pi(g)v+b(g) \tout g\in G, v\in {\mathcal H}.$$
Conversely, let $\al$  be a continuous action of  $G$ on $\H$  by affine isometries.
Denote by  $\pi(g)$ and $b(g)$ the linear part and the translation part of $\al(g)$ for  $g\in G.$
Then $\pi$ is a unitary (or orthogonal) representation of $G$ on $\H,$ $b$ is a $1$-cocycle in $Z^1(G,\pi),$
and $\al=\alpha_{\pi, b}.$ 
 For all this,  see Chapter 2 in \cite{BHV}.

The following notion of irreducibility of affine actions  was introduced in \cite{Ner} and  further  studied in \cite{BPV}.
\begin{definition}
\lb{Def-IAA}
An  affine isometric action  $\alpha$ of  $G$ on the 
complex or real Hilbert space $\H$ is {\it irreducible} if $\H$ has no non-empty, closed and  proper $\al(G)$-invariant affine subspace.
\end{definition}

First examples of irreducible affine isometric actions arise as actions $\alpha_{\pi, b},$
where $\pi$ is an irreducible unitary representation of $G$ with non trivial 
$1$-cohomology and $b\in Z^1(G,\pi)$ a cocycle which is not a coboundary.
By \cite[Theorem 0.2]{Shalom2},  such a pair $(\pi, b)$ always exists,
provided $G$ does not have Kazhdan's Property (T).
A remarkable feature of  irreducible affine isometric actions 
of a locally compact group $G$ is that they
remain irreducible under restriction to  ``most" lattices in $G$ 
(see \cite[3.6]{Ner},  \cite[Theorem 4.2]{BHV}),
whereas this is not true in general for irreducible unitary representations.

 Let $b\in Z^1(G,\pi)$.   Observe that $\Span(b(G))$ is  $\al_{\pi, b}(G)$-invariant.
So, for  $\al_{\pi, b}$ to be irreducible, it is necessary that 
$\Span(b(G))$ is dense in $\H.$ This condition is not sufficient
(see \cite[Example 2.4]{BPV}; however, see also   Proposition~\ref{Pro-IrrDens} below).
  The following corollary  of Theorem~\ref{Theo1} relates harmonic cocycles
  to this question.

  \begin{corollary}
  \lb{Theo1-Cor1} 
   Let $G, (\pi,\H),$ and $\mu$ be as in Theorem~\ref{Theo1}.
  Let $b\in Z^1(G, \pi)$ and $P_{\Har} b$ its projection on $\har (G,\pi).$
 
 \n
 (i) If  $\Span (P_{\Har} b (G))$ is dense in 
 $\H,$ then   the affine action $\alpha_{\pi, b}$ is irreducible.
 
 \n
 (ii)   Assume that $B^1(G,\pi)$ is closed; if the affine action $\alpha_{\pi, b}$ is irreducible,
 then $\Span (P_{\Har} b (G))$ is dense in   $\H.$ 
\end{corollary}
  \begin{remark}
   \lb{Rem1}
 \n
 (i)    Point (ii) in Corollary~\ref{Theo1-Cor1}  does not hold in general when $B^1(G,\pi)$ is not closed;
   indeed, let $G=\FF_2$  denote the free group on $2.$ generators. Then $H^1(G, \pi)\neq 0$
   for every unitary representation $\pi$ of $G$ (see \cite[\S 9, Example~1]{Gui2}). On the other hand,
      there exists an irreducible unitary representation $\pi$ of $G$  with  $\overline{H}^1(G, \pi)=0$
     (see \cite[Theorem 1.1]{MV}), so that  $\har(G, \pi)=0$ for any cohomologically adapted probability measure $\mu$
     on $G.$
  Now, let $b$ be a $1$-cocycle in  $Z^1(G,\pi)$ which is not a coboundary. Then   the affine action $\al_{\pi, b}$ 
  is irreducible. 
     
   \n
   (ii)  Although we will not need it, we will give an explicit formula for  the projection  $P_{\Har}: Z^1(G,\pi) \to \har(G,\mu)$  in the case where $B^1(G,\pi)$ is closed (see Proposition~\ref{Rem-ExplicitFormulaProjection} below).
 \end{remark}

  In view of Corollary~\ref{Theo1-Cor1}, it is of interest to know when $B^1(G,\pi)$ is closed.
 Write $\H= \H^G\oplus \H^0,$ where 
  $\H^G$ is the space of $\pi(G)$-invariant vectors in $\H$ and  $\H^0$ its orthogonal complement.
  Let $\pi^0$ denote the restriction of $\pi$ to $\H^0$. 
 Observe that  $B^1(G,\pi^0)=B^1(G,\pi)$ and that $Z^1(G,\pi^0)$ is closed in $Z^1(G,\pi);$
 so, the following result is both a (slight) strengthening and a consequence of Th\'eor\`eme~1 in \cite{Gui2}.
 \begin{proposition}
 \lb{Pro-Gui}
 \textbf{(\cite{Gui2})}
 Let $(\pi, \H)$ be an orthogonal or  unitary representation of the $\sigma$-compact group $G.$ 
  Then $B^1(G,\pi)$ is closed   in $Z^1(G,\pi)$ if and only if $(\pi^0, \H^0)$ does not weakly contain  
  the trivial representation $1_G.$
  \end{proposition}
  
    Our approach to the proof of Theorem~\ref{Theo1}  uses the fact, observed in \cite[\S 3.1]{BPV}
    and  \cite{BV} that
$\overline{H}^1(G,\pi)$ is, in a natural way, a module over the (real or complex) von Neumann algebra $\pi(G)',$ 
which is the commutant of $\pi(G)$ in $\B(\H);$ see Section~\ref{S:Cocycle-VN}.
Viewing, as we do, $\overline{H}^1(G,\pi)$ as the Hilbert space $\har (G,\pi),$ one is lead to the study
of $\har (G,\pi)$ as a Hilbert module over  $\pi(G)'.$ 

For instance, if  $\M:=\pi(G)'$ is a finite von Neumann algebra
(that is,  if there exists a faithful finite trace on $\M$)  then, we can define (as in \cite[Definition p.138]{GHJ} or \cite[p. 327]{Bek})  the \emph{von Neumann dimension}   of $\overline{H}^1(G,\pi) $ as   
$$\dim_\M\overline{H}^1(G,\pi):= \dim_\M \har (G,\pi) \in[0,+\infty)\cup \{+\infty\};$$
for more details, see Section~\ref{S:Cocycle-VN}.
It is worth mentioning that in case $\pi$ is the regular representation of a discrete group $G,$
$\dim_\M \overline{H}^1(G,\pi)$ coincides with  $\beta_2^1(G),$ the $L^2$-Betti number of $G$
(see \cite[Proposition~2]{BV}).

We now give some applications  of von Neumann techniques to the problem of the existence of an irreducible  affine isometric action of $G$ with a given linear part $\pi.$ First, using Corollary~\ref{Theo1-Cor1},
 we can reformulate Corollary~3.7 from \cite{BPV} in our setting.
Recall that a vector $v$ in a Hilbert module  over a von Neumann algebra $\M$ is a  \emph{separating vector} for 
$\M$ if $Tv= 0$ for $T\in \M$ implies $T=0.$

\begin{proposition}
\lb{Pro-Totalisateur}
\textbf{(\cite{BPV})}

\n
(i) Assume that $\M=\pi(G)'$ has a separating vector $b$ in $\har(G,\pi).$ Then $\alpha_{\pi, b}$ is irreducible.

\n
(ii) Assume $B^1(G,\pi)$ is closed and that $\alpha_{\pi, b}$ is irreducible for some $b\in \har(G,\pi).$ Then
$b$ is  a separating vector for $\M.$
\end{proposition}

For an application of the previous criterion in  the case where $G$ is a discrete finitely generated group
and $\pi$ a subrepresentation of a multiple of the regular representation of $G$, see \cite[Theorem~4.25]{BPV}.
We extend this result to arbitrary factor representations, that is, to  unitary  representations 
$(\pi, \H)$ such that the von Neumann subalgebra $\pi(G)''$ of  $\B(\H)$ generated by $\pi(G)$ 
is a factor (equivalently, such that $\pi(G)'$ is a factor).
Concerning general facts about factors, such as their type classification, see \cite{Dixmier-VN}.

\begin{theorem}
\lb{Theo2}
Let $(\pi, \H)$ be a factor representation of the compactly generated
locally compact group $G$ on the separable complex Hilbert space $\H.$
 Assume that $B^1(G,\pi)$ is closed in $Z^1(G, \pi).$ 
Set $\M:=\pi(G)'$ and let  $\mu$ be a cohomologically adapted probability measure on $G.$ 
Depending on the type of $\M,$ there exists $b\in Z^1(G, \pi)$ such that $\alpha_{\pi, b}$ is irreducible
if and only  if: 

\begin{itemize}
\item[(i)] the factor $\M$ is   of type $I_\infty$ or of type $II_\infty$
and  its commutant   in $\B(\har(G, \pi))$ is  of infinite type (that is, of type $I_\infty$ or $II_\infty,$ respectively);
\item[(ii)] the factor  $\M$ is of finite type (that is,  of   type $I_n$ for $n\in \NN$ 
or of type $II_1$) and $\dim_\M \har(G, \pi) \geq 1;$
\item[(iii)] the factor $\M$ is   of type $III$ and $\har(G, \pi)\neq \{0\}.$ 
\end{itemize}
\end{theorem}

\begin{remark}
\lb{Rem-AA-VNA}
 Let $(\pi, \H)$ be a  unitary representation of $G$ such that 
$B^1(G,\pi)$ is closed in $Z^1(G,\pi)$;  let 
$$\pi=\int^{\oplus}_\Omega\pi_\omega d\nu(\omega)$$
 be the central  integral decomposition 
of $\pi$, so that the $\pi_\omega$'s  are mutually disjoint factor representations of $G$
(see \cite[Theorem~8.4.2]{Dixmier-C*}). 
One checks that one has a corresponding decomposition of $\har (G,\pi)$ as a direct integral  of
Hilbert spaces:
$$
\har (G,\pi)=\int^{\oplus}_\Omega\har (G,\pi_\omega) d\nu(\omega).
$$
Moreover, $B^1(G,\pi_\omega)$ is closed in $Z^1(G,\pi_\omega)$
and there exists a separating vector for $\pi(G)'$ in $\har (G,\pi)$ if and only if 
 there exists a separating vector for $\pi_\omega(G)'$ in $\har (G,\pi_\omega)$
 for $\nu$-almost every $\omega.$
 So, Theorem~\ref{Theo2} can be used to  check  
 the existence of an irreducible affine with \emph{any} unitary representation $\pi$ as linear part
 (provided $B^1(G,\pi)$ is closed in $Z^1(G,\pi)$).
\end{remark} 

As an illustration of the use of Theorem~\ref{Theo2}, we will 
treat the example  of a wreath product of the form  $\Ga=G\wr \ZZ$
and a unitary representation $\pi$ of $\Ga$ which factorizes
 through a representation of $G;$ the reduced cohomology of such groups
 was considered in \cite[\S 5.4]{Shalom}.

 \begin{theorem}
 \label{Theo3}
 Let $G$ be a finitely generated group, 
 and let  $(\pi, \H)$ be a unitary representation of the wreath product $\Ga=G\wr \ZZ$
 in the separable Hilbert space $\H.$ 
 Assume that $\pi$ factorizes  through $G$ and that $H^1(G,\pi)=0$.

 \n
 (i)  For a suitable cohomologically adapted probability measure 
 $\mu$ on $\Ga,$ the space $\har (\Ga, \mu)$ can be identified,
 as a module over $\pi(\Ga)'$, with the Hilbert space $\H.$
 
 \n
(ii) There exists  an irreducible affine action of $\Ga$ with linear part $\pi$
 if and only if the representation $(\pi, \H)$ is cyclic.
 
\n
(iii)  Assume that  $G$ is not virtually abelian (that is, $G$  does not have
 an abelian normal subgroup of finite  index). Then  $G$  has a factorial representation $\pi$ 
 for which  $\pi(G)'$ is of any possible type.
 
 \end{theorem}
\begin{remark}
 \lb{Rem-WR}
 \n
 (i) When  $\pi$ is a factor representation, a necessary and  sufficient condition for the existence
 of a cyclic vector for $\pi(G)$  (equivalently, a separating vector for $\pi(G)'$) in $\H$ is given in 
 Theorem~\ref{Theo2}, with $\H$ replacing $\har(G,\mu)$ there.
 
 \n
 (ii)  By the Delorme-Guichardet theorem (\cite[Theorem 2.12.4]{BHV}), the condition $H^1(G,\pi)=0$ 
is satisfied for every unitary representation $\pi$ of $G$ if (and only if) $G$ has Kazhdan's property (T).  
\end{remark}

  \section{The space of harmonic cocycles as  a von Neumann algebra module}
   \label{S:Cocycle-VN}
 
  Let $G$ be a locally compact group  which is generated by  a compact  subset $Q$, which we assume to be a 
  symmetric neighbourhood of the identity $e\in G$. Let $(\pi, \H)$ be an orthogonal or unitary representation of $G.$
     The map 
     $$b\mapsto \Vert b\Vert_Q=\sup_{x\in Q} \Vert b(x)\Vert$$
 is a norm  which generates the topology of uniform convergence on compact subsets and  for which  $ Z^1(G,\pi)$ is a Banach space.
   
    Let $\M:=\pi(G)'$ be  the commutant of $\pi(G)$ in $\B(\H),$ that is, 
$$\M=\{T\in{{B}({\H})}: T\pi(g)=\pi(g)T \quad \text{for all} \quad  g \in G\};$$ 
 $\M$ is a  (real or complex) von Neumann algebra, that is,
 $\M$ is a unital self-adjoint subalgebra of $\B(\H)$ which is closed for the weak (or strong) operator topology.

As observed in \cite[\S 3.1]{BPV}), $H^1(G,\pi)$  is a module over $\M$;  indeed, if $b \in Z^1(G,\pi)$ and $T\in\pi(G)'$,  then $Tb \in Z^1(G,\pi),$   where $Tb$ is defined by 
$$Tb(g)=T(b(g)) \tout g\in G;$$
 moreover, $T\partial_v= \partial_{Tv}$  for every vector $v\in \H.$
 
 Let $\mu$  be a cohomologically adapted probability measure on $G$ (Definition~\ref{Def-AdmissibleMeasure}).
  We consider the Hilbert space $ L^2(G,\H, \mu)$ of measurable mappings  $F: G\to \H$ such that 
  $$\Vert F\Vert_2^2:= \int_G \Vert F(x)\Vert^2 d\mu(x) <\infty.$$
  Then every $b\in Z^1(G,\pi)$ belongs to $L^2(G,  \H, \mu);$  indeed,  the cocycle relation shows that 
  $$\Vert b(x)\Vert \leq |x|_Q \Vert b\Vert_Q \tout x\in G,$$
 and hence
 $$
    \Vert b\Vert _2^2  \leq \Vert b\Vert_Q^2 \int_{G} |x|_Q^2 d\mu(x) <\infty.
 $$
    
 In fact, the norms $\Vert \cdot \Vert_2$ and $\Vert \cdot \Vert_Q$ on $Z^1(G,\pi)$ are equivalent 
 (see \cite[Lemma 2.1]{Anna-Ozawa}).
  So, we can (and will) identify $Z^1(G,\pi)$ with a closed
  subspace of the Hilbert space $L^2(G,\H, \mu).$ 
  
  The  von Neumann algebra $\M$ acts on $\H$ in the tautological way and 
   on $L^2(G,\H, \mu)$ by
    $$
  TF(g)= T(F(g)) \tout T\in \pi(G)', F\in L^2(G,\H, \mu), g\in G,
  $$
  preserving $Z^1(G,\pi)$ and $B^1(G,\pi).$
  Since the operator $M_\mu: Z^1(G,\mu)\to \H$ is equivariant for these actions,
  $\har(G,\pi)=\ker M_ \mu$ as well as its orthogonal complement $\overline{B^1(G,\pi)}$  are 
  modules over $\M.$ 
  
  The image of $\M$ in $\B(L^2(G,\H, \mu))=\B(L^2(G,\mu))\otimes \H$ is 
  $$\widetilde{\M}=  I \otimes \pi(G)',$$
  which is a von Neumann algebra isomorphic to $\M.$
  The orthogonal projection $P_\Har: L^2(G,\H, \mu)\to \har(G,\pi)$ belongs
   to the commutant  
  $$\widetilde{\M}'=\B(L ^2(G,\mu))\otimes \pi(G)''$$
  of $\M$ in $\B(L^2(G,\H, \mu)),$ where $\pi(G)''$ is the von Neumann algebra generated by $\pi(G)$
  in $\B(\H).$ The commutant of $\M$ in $\har (G,\pi))$ is then the reduced von Neumann algebra 
  (see Chap.1, \S, Proposition~1 in \cite{Dixmier-VN}) 
  $$
  P_{\Har} \widetilde{\M}' P_{\Har}= P_{\Har}(\B(L ^2(G,\mu))\otimes \pi(G)'')P_{\Har}.
  $$
  
 Assume now that $\M$ is a finite von Neumann algebra, with  
  faithful  normalized  trace $\tau.$ 
Let $L^2(\M)$ be the Hilbert space  obtained from $\tau$ by 
the GNS construction. We  identify $\M$  with the subalgebra of $\B(L^2({\M}))$ of operators
given by left multiplication with elements from $\M.$ The commutant of $\M$ in
$\B(L^2(\M))$  is $\M '=J\M J,$ where $J : L^2(\M) \to L^2(\M)$ is the conjugate linear isometry
which extends the map $\M\to \M, x\mapsto x^*$. The trace on  $\M'$, again denoted by $\tau,$
is defined by  $\tau(JxJ)= \tau (x)$ for $x\in \M.$

The $\M$-module $ L^2(G,\H, \mu)$ can be identified with an $\M$-submodule of
$L^2(\M)\otimes \ell^2$, with $\M$ acting on  $L^2(\M)\otimes \ell^2$  by  
$T\mapsto T\otimes I.$
The  orthogonal projection $Q: L^2(\M)\otimes \ell^2 \to L^2(G,\H, \mu)$ 
belongs to the commutant of $\M$
in ${\B}(L^2({\M})\otimes \ell^2)$, which is  ${\M}' \otimes {\B}(\ell^2).$
The projection $P= P_\Har \circ Q$ belongs therefore to  the commutant of $\M$ in ${\B}(L^2({\M})\otimes \ell^2).$

Let $\{e_n\}_n$ be a basis of $\ell^2$ and  let $(P_{ij})_{i,j}$ be the matrix of 
$P$ with respect to the decomposition  $L^2({\M})\otimes  \ell^2= \oplus_{i} (L^2({\M}) \otimes e_i).$
Then every $P_{ij}$ belongs to ${\M} '$ and 
the von Neumann dimension  of the $\M$-module $\har(G,\pi)$
is
$$
 {\rm dim}_{\M}{\H} = \sum_{i} \tau(P_{ii}).
  $$

 \section{Proofs of the main results}
  \subsection{Proof of Theorem~\ref{Theo1}}
   \lb{SS:ProofTheo1-Cor}
  
  Let $b_0\in \har(G,\pi).$    Let $b_1\in \overline{B^1(G,\pi)}$ and set $b:=b_0+b_1$.
 We claim that  $b_0(G)$
  is contained in the closure of $\Span (b(G)).$
  
  Indeed, let $\K$ denote the closure of $\Span (b(G))$
  and $P_{\K}: \H\to \K$ the corresponding orthogonal projection.
  Since $\K$ is $\pi(G)$-invariant,  $P_{\K}$ belongs to the commutant  $\pi(G)'$
  of $\pi(G).$ Therefore (see Section~\ref{S:Cocycle-VN}), $\P_{\K}b_0$ is contained in   $\har(G,\pi)$ and 
   $\P_{\K}b_1$ is contained in $\overline{B^1(G,\pi)}.$ On the other hand, since
    $b$ take its values in $\K,$ we have that 
$$P_{\K} b= b= b_0+b_1.$$  It follows that 
  $\P_{\K}b_0 = b_0$ and $\P_{\K}b_1= b_1$.
  Therefore, 
  $$b_0(G) \subset \K= \overline{\Span(b(G))},$$
   as claimed. $\bsq$

 \subsection{A characterization of  irreducible  affine isometric actions }
   \lb{SS:CharIAA}
  We will need for the proof of Corollary~\ref{Theo1-Cor1}   one of the several characterizations of irreducible  affine actions from Proposition~2.1 in \cite{BPV};  for the convenience of the reader, we give   a direct and short argument.
  
  \begin{proposition}
  \lb{Pro-IrrDens}
  \textbf{(\cite{BPV})}
 For  $b\in Z^1(G,\pi),$ the following properties are equivalent:
 
\n
(i) the action $\al=\al_{\pi, b}$ is irreducible;

\n
(ii) the linear span of  $(b+\partial_v)(G))$ is dense in $\H$ for every $v\in \H.$
\end{proposition}

\proof
Observe that 
$$\al_{\pi, b+\partial_v}(g)=t_{-v}\circ \al_{\pi, b}(g)\circ t_v \tout g\in G, \  v\in \H,$$
 where $t_v$ is the translation by $v.$ So, $\al_{\pi, b}$
is irreducible if and only if $\al_{\pi, b+\partial_v}$ is irreducible.
This shows that (i) implies (ii). 

To show the converse implication, let $F$ be a non empty closed $\alpha_{\pi, b}(G)$-invariant
affine  subspace of $\H.$ Then $F= v + \K$ for a vector $v\in \H$
and a   closed linear subspace $\K$ of $\H.$  
  Set $b_0:= b+ \partial_{v}.$ 
  Then
  $$v+ b_0(g)=\alpha_{\pi, b}(g)v\in F \tout g\in G,$$
 and  $b_0(G)$ is hence contained in  $\K.$  Therefore, $\K=\H, $ since $\Span (b_0(G))$
 is dense in $\H.\bsq$

 \subsection{Proof of Corollary~\ref{Theo1-Cor1}} 
 Let $b\in Z^1(G,\pi)$ and set 
$b_0:=P_{\Har} b \in \har(G,\pi).$ 
 
\n
(i)   Assume that
 $\Span (b_0(G))$  is dense in $\H$.
  By Theorem~\ref{Theo1},   the  linear span of $(b+\partial_{v})(G)$ is dense for every $v\in \H,$
  and Proposition~\ref{Pro-IrrDens} shows
  that  $\alpha_{\pi, b}$ is irreducible.

  \n
  (ii)  Assume now that $B^1(G,\pi)$ is  closed in $Z^1(G,\pi)$ and that 
  $\alpha_{\pi, b}$ is irreducible. 
  Write $b=b_0+\partial_{v_0}$ for
  $b_0=P_{\Har} b$ and $v_0\in \H.$  Then $\alpha_{\pi, b_0}=\alpha_{\pi, b-\partial_{v_0}}$
  is also irreducible, by Proposition~~\ref{Pro-IrrDens}; hence, $\Span(b_0(G))$ is dense.$\bsq$

     \subsection{Proof of Theorem~\ref{Theo2}}
  \lb{SS:Theo2}
  Let $(\pi, \H)$ be a unitary representation of $G;$ 
  we assume that $B^1(G,\pi)$ is closed in $Z^1(G, \pi).$ 
Let  $\mu$ be a cohomologically adapted probability measure on $G.$ 

In view of Proposition~\ref{Pro-Totalisateur}, we have to investigate under which conditions
$\M=\pi(G)'$ has a separating vector in $\har(G,\pi).$ We may assume that $\har(G,\pi)\neq \{0\}.$

Observe that a vector  in $\har(G,\pi)$  is separating for $\M$  if and only if 
it is cyclic for the commutant $\N$ of $\M$ in $ \B(\har(G,\pi)).$  
Three  cases cases can occur.

\n
$\bullet$ \emph{First case:} $\N$  is an infinite factor. Then 
  $\M$ always has a separating vector  (see Corollaire 11 in Chap. III, \S 8 of \cite{Dixmier-VN}).

\n
$\bullet$ \emph{Second case:} $\N$  is a finite factor and $\M$ is an infinite factor. 
Then $\N$ has a cyclic vector in $\har(G,\pi)$ 
if and only if $\dim_\N \har(G,\pi) \leq 1$ (see  \cite[Corollary 1]{Bek}).
For  this to happens a necessary condition is that $\M$ is a finite factor. So, $\M$
has no separating vector.

\n
$\bullet$ \emph{Third case:} $\N$  and  $\M$  are  finite factors. 
In this case, we have  (see \cite[Prop. 3.2.5]{GHJ})
$$\dim_\M \har(G,\pi)\dim_\N\har(G,\pi)=1;$$
 hence, $\M$ has a separating vector in $\har(G,\pi)$  if and only if 
 $$\dim_\M \har(G,\pi)\geq 1.$$
 Claims (i), (ii), and (iii) follow from this discussion.$\bsq$

\subsection{Proof of Theorem~\ref{Theo3}}
\lb{S:Examples}
 We  first  consider the general  case of the wreath product  $\Ga=G\wr H$  of two finitely generated  groups $G$ and $H.$ Recall
 that  $\Ga= G\ltimes H^{(G)},$
 for  $H^{(G)}=\bigoplus_{g\in G} H$ and $G$ acts on  $H^{(G)}$ by shifting the copies of $H.$
 We view $H$ as a subgroup of $\Ga$, by identifying it with the copy of $H$ inside $H^{(G)}$ indexed by $e.$
 
Let $S_1$ and $S_2$ finite symmetric generating sets for $G$ and $H,$ respectively.
Then $S_1\cup S_2$ is a finite symmetric generating set for $\Ga.$
Let $\mu_1$ and  $\mu_2$  be  cohomologically adapted probability measures on  $G$ 
 and  $H$ respectively. Then $\mu=\dfrac{1}{2}(\mu_1+\mu_2)$ is 
a cohomologically adapted probability measure on $\Ga.$

 Let $(\pi, \H)$ be a unitary representation of $G,$
viewed as representation of $\Ga.$
We have   orthogonal $\pi(\Ga)$-invariant decompositions
$$ 
\ell^2(\Ga, \H, \mu)= \ell^2(G,\H, \mu_1) \oplus  \ell^2(H,\H, \mu_2)
$$
and 
$$ 
\har (\Ga, \pi)=  \Har_{\mu_1} (G, \pi)  \oplus  \Har_{\mu_2} (H, \pi).
$$
Since $\pi$ is trivial on $H,$ the space  $Z^1(H, \pi)$ coincides with  
the set $\Hom(H, \H)$ of homomorphisms $H\to \H.$
Observe that every  $b\in  \Hom (H, \H)$ is $\mu_2$-harmonic, 
 since
 $$
 \sum_{h\in H } b(h)\mu_2(h)=  \sum_{h\in H } b(-h)\mu_2(h)= -\sum_{h\in H } b(h)\mu_2(h).
 $$
 Hence, $\Har_{\mu_2} (H, \pi)= \Hom (H, \H)$ (alternatively, this follows
 from the fact that $B^1(H, \pi)=B^1(H, 1_H)$ is trivial); therefore, we have
$$ 
\har (\Ga, \pi)=  \Har_{\mu_1} (G, \pi)  \oplus  \Hom (H, \H).
$$

 We specialize  by taking $H=\ZZ;$ then $\Hom(H, \H)$ can be identified with $\H$
 and we have
 $$ 
\har (\Ga, \pi)=  \Har_{\mu_1} (G, \pi)  \oplus  \H;
$$
 moreover, the action of the von Neumann algebra $\pi(\Ga)'=\pi(G)'$ on $\har (G, \mu)$
 corresponds to the direct sum of the actions of $\pi(G)'$ on $\Har_{\mu_1} (G, \mu_1)$ and 
 on  $\H$.

In particular, when the $1$-cohomology $H^1(G,\pi)$ is trivial,
we have $$\har (\Ga, \pi)=\H,$$  so that Claim  (i)  is proved.
Claim (ii) follows from  Proposition~\ref{Pro-Totalisateur}.

  To show Claim (iii), assume that  $G$ is not virtually abelian. Then $G$ is not of type $I,$ by Thoma's theorem (\cite[Satz~6]{Thoma}).
  
 First,  observe that $G$ has an irreducible  unitary representation $\sigma$ of infinite dimension;
 indeed, otherwise, $G$ would be a liminal (or CCR) group  and hence  of type $I$, by \cite[13.9.7]{Dixmier-C*}.
 Set  $\pi=n\sigma$, a multiple of $\sigma$ for $n\in \NN$ or $n=\infty;$
 then $\pi(G)'$ is of type $I_n.$
 
 Next, since $G$ is not of type $I,$
 $G$ has a factorial representation $\pi$ such that both $\pi(G)''$ and  $\pi(G)'$ are of type $II_1$,
 by \cite[Lemma 19]{Thoma}. Then $\rho:=\infty \pi$   is factorial and $\rho(G)'$ is of type $II_\infty.$
 
 Finally, $G$ has a factor representation such that $\pi(G)''$ (and hence $\pi(G)'$)  is of type $III$,
 by Glimm's theorem \cite[Theorem~1]{Glimm}). $\bsq$
 
     \section{An explicit formula for the projection on harmonic cocycles}
     \lb{SS:ExplicitFormula}
  We give an explicit formula for the orthogonal projection  $P_{\Har}$
  in terms of an averaging (or Markov) operator associated to $\mu,$
  in the case where $B^1(G,\pi)$ is  closed.  
  
  Consider the  operator $\pi^0(\mu) \in \B(\H^0)$ defined by 
  $$\pi^0(\mu)v= \int_G \pi(x)v d\mu(x) \tout v\in \H^0.$$
  The operator $\pi^0(\mu)-I: \H^0\to \H^0$ is invertible  if and only 
  if $\pi^0$ does not  weakly contain  
  the trivial representation $1_G$  (see Proposition G.4.2 in \cite{BHV}); 
   in view of Proposition~\ref{Pro-Gui}, this is the case if and only if $B^1(G,\pi)$ is closed.
  \begin{proposition}
  \lb{Rem-ExplicitFormulaProjection}
  Assume that $B^1(G,\pi)$ is closed. 
 For $b\in Z^1(G,\pi),$ we have $P_{\Har} b= b- \partial_v,$ where
  $$
 v= (\pi^0(\mu)-I)^{-1} (M_\mu(b)).
  $$
  \end{proposition}
  \proof
  Indeed, observe first that $M_\mu(b)\in \H^0;$ indeed, for every $w\in \H^G,$ we have 
  \begin{align*}
  \langle M_\mu(b), w\rangle &= \int_G \langle b(x), w\rangle d\mu(x)=  \int_G \langle b(x), \pi(x) w\rangle d\mu(x)\\
  &=\int_G \langle \pi(x^{-1})b(x), w\rangle d\mu(x)= -\int_G \langle b(x^{-1}), w\rangle d\mu(x)\\
  &=- \int_G \langle b(x), w\rangle d\mu(x)= - \langle M_\mu(b), w\rangle.
 \end{align*}
  Moreover, for $v= (\pi^0(\mu)-I)^{-1} (M_\mu(b)),$ we have
  \begin{align*}
  M_\mu( \partial_v)&=\int_G (\pi(x)v-v) d\mu(x)=  (\pi^0(\mu)-I)v= M_\mu(b).\bsq
 \end{align*}

 \end{document}